\documentclass[11pt,a4paper]{article}
\usepackage{amsmath,amssymb}
\usepackage{bm}
\usepackage{ascmac}
\usepackage[dvipdfmx]{graphicx}
\usepackage{enumerate}
\usepackage{tikz}
\usepackage[left=3.2cm, right=3.2cm]{geometry}

\newcommand{\Rone}{\mathrm{I}}
\newcommand{\Rtwo}{\mathrm{II}}
\newcounter{Cntr}
\usepackage{amsthm}
\newtheorem{thm}{Theorem}[section]

\newtheorem{lem}[thm]{Lemma}
\newtheorem{cnj}[thm]{Conjecture}

\newtheorem{cla}[thm]{Claim}

\theoremstyle{definition}
\newtheorem{defn}[thm]{Definition}
\newtheorem*{prf}{Proof}

\newtheorem{case}{Case}
\newtheorem*{subcase1}{Case 2.1}
\newtheorem*{subcase2}{Case 2.2}

\usepackage{latexsym}
\def\qed{\hfill $\Box$}

\definecolor{mygray}{gray}{0.9}

\begin{document}

\title{Forbidden subgraphs and 2-factors in 3/2-tough graphs}

\author{%
    Masahiro Sanka\thanks{E-mail address: \texttt{sankamasa@keio.jp}}\\
    Department of Mathematics, Keio University, \\
    3-14-1 Hiyoshi, Kohoku-ku, Yokohama 223-8522, Japan
} 

\date{}

\maketitle

\begin{abstract}
    A graph $G$ is $H$-free if it has no induced subgraph isomorphic to $H$, where $H$ is a graph. 
    In this paper, we show that every $\frac{3}{2}$-tough $(P_4 \cup P_{10})$-free graph has a 2-factor.
    The toughness condition of this result is sharp. Moreover, for any $\varepsilon>0$ there exists a $(2-\varepsilon)$-tough $2P_5$-free graph without a 2-factor.
    This implies that the graph $P_4 \cup P_{10}$ is best possible for a forbidden subgraph in a sense.
\end{abstract}


\section{Introduction}\label{sec:intro}

We begin with a few definitions.
The terminology not defined here can be found in \cite{Diestel-2017}, and additional definitions will be given as needed.

All graphs in this paper are finite, undirected and simple.
Let $G$ be a graph with vertex set $V(G)$ and edge set $E(G)$.
For $x \in V(G)$, we denote by $N_G(x)$ and $d_G(x)$ the neighborhood and the degree of $x$ in $G$, respectively.
We denote by $\delta(G)$ and $\omega(G)$ the minimum degree and the number of components of $G$, respectively.

Let $U \subset V(G)$.
Then we define that $N_G(U)=\left( \bigcup_{x \in U}N_G(x) \right) \setminus U$.
If any two vertices in $U$ are adjacent in $G$, we call $U$ a {\it clique}.
Also, if any two vertices in $U$ are not adjacent in $G$, we call $U$ an {\it independent set}. 
Let $H$ be a subgraph of $G$.
Then we write $N_G(H)$ for $N_G(V(H))$.
If $E(H)=\{uv \mid u,v \in V(H) \text{ and }uv \in E(G)\}$, we say that $H$ is an {\it induced subgraph} of $G$, and $V(H)$ induces $H$.
The subgraph of $G$ induced on $V(G) \setminus U$ is denoted by $G-U$.
For notational simplicity we write $G-x$ for $G-\{x\}$.

Let $X,Y \subset V(G)$, and $x \in V(G)$. 
We use $e_G(x,Y)$ to denote the number of edges of $G$ joining $x$ to a vertex in $Y$, $e_G(X,Y)$ to denote $\sum_{x \in X}e_G(x,Y)$, $N_Y(x)$ to denote $Y \cap N_G(x)$ and $N_Y(X)$ to denote $Y \cap N_G(X)$.
If $e_G(x,Y) \geq 1$ $(e_G(X,Y) \geq 1)$, we say that $x$ and $Y$ ($X$ and $Y$) are {\it adjacent}.

For two graphs $G$ and $H$ and a positive integer $m$, we define that $G \cup H$ is the disjoint union of $G$ and $H$ and $mG$ is the disjoint union of $m$ copies of $G$.
We use $K_n$ and $P_n$ to denote the complete graph and the path with $n$ vertices, respectively.

A {\it $k$-factor} of a graph $G$ is a spanning subgraph of $G$ in which every vertex has degree $k$.
In particular, a 2-factor is a spanning subgraph in which every component is a cycle. 
We call a connected 2-factor a {\it hamiltonian cycle}, and say that $G$ is {\it hamiltonian} if $G$ has a hamiltonian cycle.

The {\it toughness} of a graph $G$, denoted by $t(G)$, is defined by
\[t(G)=\min \left\{\frac{|U|}{\omega(G-U)} \mid U \subset V(G) \text{ and } \omega(G-U) \geq 2\right\}\]
or $t(G)=\infty$ if $G$ is complete.
For a real number $t$, if $t(G) \geq t$, then we say that $G$ is a {\it $t$-tough} graph.
Clearly, a graph $G$ is $t$-tough if and only if $t \cdot \omega(G-U) \leq |U|$ for any subset $U \subset V(G)$ with $\omega(G-U) \geq 2$.

The notion of toughness was introduced in the study of hamiltonian graphs by Chv\'{a}tal \cite{Chvatal-1973}. 
It is well known that every hamiltonian graph is 1-tough, but the converse does not hold. 
Chv\'{a}tal has conjectured that the converse holds at least in an approximate sense.

\begin{cnj}[\cite{Chvatal-1973}, 1973]\label{cnj:Chvatal}
    There exists a constant $t_0$ such that every $t_0$-tough graph on at least three vertices is hamiltonian.
\end{cnj}

In \cite{Bauer_et_al-2000}, Bauer, Broersma and Veldman constructed a $(\frac{9}{4}-\varepsilon)$-tough graph which is not hamiltonian for any $\varepsilon>0$.
Thus, the constant $t_0$ is at least $\frac{9}{4}$ if Conjecture \ref{cnj:Chvatal} is true. 
However, it is still open whether every $\frac{9}{4}$-tough graph on at least three vertices is hamiltonian.
For details of known results on Conjecture \ref{cnj:Chvatal}, we refer the reader to the survey \cite{Bauer_et_al-2006}.

Partial results related to Conjecture \ref{cnj:Chvatal} have been obtained in some restricted classes of graphs.
In particular, our work was inspired by the result on the class of split graphs.
A graph $G$ is a {\it split graph} if the vertex set of $G$ can be partitioned into a clique and an independent set. 
Kratsch, Lehel and M\"{u}ller showed the following theorem.

\begin{thm}[\cite{Kratsch_et_al-1996}, 1996]\label{thm:split graph}
    Every $\frac{3}{2}$-tough split graph on at least three vertices is hamiltonian.
    Moreover, for any $\varepsilon>0$ there exists a $(\frac{3}{2}-\varepsilon)$-tough split graph which has no $2$-factor.
\end{thm}

It has been shown that Conjecture \ref{cnj:Chvatal} is true for some superclass of split graphs, for example, spider graphs~\cite{Kaiser_et_al-2007}, chordal graphs~\cite{Chen_et_al-1998,Kabela_et_al-2017}, $2K_2$-free graphs~\cite{Broersma_et_al-2014,Shan-2019,Ota_et_al-2021}, and $(P_2 \cup P_3)$-free graphs~\cite{Shan-2021}.
However, some of the above results are not known to be best about the toughness condition, which cannot be smaller than $\frac{3}{2}$ by Theorem \ref{thm:split graph}.

In order to estimate the condition for the existence of hamiltonian cycles, it is natural to study 2-factors.
In the general case, the following theorem was proved by Enomoto, Jackson, Katerinis and Saito.

\begin{thm}[\cite{Enomoto_et_al-1985}, 1985]\label{thm:Enomoto et al}
    Let $k \geq 2$ and $n \geq k+1$ be integers with $kn$ even. 
    Then every $k$-tough graph on $n$ vertices has a $k$-factor.
    Moreover, for any $\varepsilon>0$ there exists a $(k-\varepsilon)$-tough graph on $n$ vertices with $n \geq k + 1$ and $kn$ even which has no k-factor.
\end{thm}

Theorem \ref{thm:Enomoto et al} implies that every 2-tough graph on at least three vertices has a 2-factor, and that there exists an infinite sequence of graphs without 2-factor having toughness approaching 2.
On the other hand, the toughness boundary for split graphs to have a 2-factor is $\frac{3}{2}$. 
This result was extended to the class of chordal graphs in \cite{Bauer_et_al-2000-Chordal}.

For a graph $H$, a graph $G$ {\it $H$-free} if it contains no induced subgraph isomorphic to $H$.
Note that for two graphs $H$ and $K$, if $H$ contains an induced subgraph isomorphic to $K$, then every $K$-free graph is an $H$-free graph.
So in this paper, we study as large a graph $H$ as possible, such that every $\frac{3}{2}$-tough $H$-free graph has a 2-factor, and show the following result.

\begin{thm}\label{thm:main1}
    Every $\frac{3}{2}$-tough $(P_4 \cup P_{10})$-free graph on at least three vertices has a $2$-factor.
\end{thm}

Since the class of $(P_4\cup P_{10})$-free graphs is a superclass of split graphs, the toughness condition of Theorem \ref{thm:main1} is sharp.
Furthermore, in \cite{Bauer_et_al-1994}, Bauer and Schmeichel constructed a $(2-\varepsilon)$-tough $2P_5$-free graph without 2-factor for any $\varepsilon>0$.
Since the class of $2P_5$-free graphs is a common subclass of $(P_4 \cup P_{11})$-free graphs and $(P_5 \cup P_{10})$-free graphs, the length of the forbidden paths in the assumptions of Theorem \ref{thm:main1} is best. 

The proof of Theorem \ref{thm:main1} is given in the subsequent sections.
In Section \ref{sec:pre}, we recall a well-known theorem of Tutte~\cite{Tutte-1952} for the existence of 2-factors and introduce related lemmas.
In Section \ref{sec:mainproof1}, we prove Theorem \ref{thm:main1}.
We believe that our work might be useful in proving Conjecture \ref{cnj:Chvatal} for more classes of graphs.

\section{Tutte's 2-factor theorem and related lemmas}\label{sec:pre}

In this section, we write a criterion given in \cite{Tutte-1952} by Tutte for a graph to have a 2-factor. 
Let $S$ and $T$ be disjoint subsets of $V(G)$, and $C$ be a component of $G-(S \cup T)$.
Then $C$ is said to be an {\it odd component} if $e_G(C,T) \equiv 1 \mod 2$. 
Let $h(S,T)$ be the number of odd components of $G-(S \cup T)$. 
Also, let
\[\eta(S,T)=2|S|-2|T|+\sum_{x \in T}d_{G-S}(x)-h(S,T).\]
The following theorem is the case $k=2$ of Tutte's $k$-factor theorem.

\begin{thm}[\cite{Tutte-1952}, 1952]\label{thm:2-factor}
    For a graph $G$, each of the following holds:
    \begin{enumerate}[$(i)$]
        \item For every disjoint sets $S,T \subset V(G)$, $\eta(S,T)$ is even.
        \item $G$ has a $2$-factor if and only if $\eta(S,T) \geq 0$ for every disjoint sets $S,T \subset V(G)$.
    \end{enumerate}
\end{thm}

Let $G$ be a graph without 2-factor. 
Then by Theorem \ref{thm:2-factor}, there exists a pair $(S,T)$ of disjoint subsets of $V(G)$ with $\eta(S,T) \leq -2$, and we call such a pair a {\it Tutte pair} for $G$.  
We call a Tutte pair $(S,T)$ a {\it special Tutte pair} if among all the Tutte pairs for $G$, 
\begin{enumerate}[(a)]
    \item $|S|$ is maximum; 
    \item $|T|$ is minimum subject to (a); and 
    \item $h(S,T)$ is minimum subject to (b).
\end{enumerate}
We use the following lemma for a special Tutte pair.
This lemma was proved, for example, in \cite{Kanno_et_al-2019}.

\begin{lem}[\cite{Kanno_et_al-2019}, Lemma 3.2, Lemma 3.3]\label{lem:sp Tutte pair1}
    Let $G$ be a graph without $2$-factor, and let $(S,T)$ be a special Tutte pair of $G$. 
    Then each of the followings holds:
    \begin{enumerate}[$(i)$]
        \item The set $T$ is independent in $G$.
        \item Each vertex $x \in T$ is adjacent to exactly $d_{G-S}(x)$ odd components of $G-(S \cup T)$.
        \item If $C$ is an odd component of $G-(S \cup T)$, then for any $y \in V(C)$, $e_G(y,T) \leq 1$.
        \item Let $x \in T$ be a vertex with $d_{G-S}(x) \geq 2$.
        If $x$ is adjacent to an odd component $C$ of $G-(S \cup T)$, then $|V(C)| \geq 3$.
    \end{enumerate}
\end{lem}

\section{Proof of Theorem \ref{thm:main1}}\label{sec:mainproof1}

In order to prove Theorem \ref{thm:main1}, suppose that $G$ is a $\frac{3}{2}$-tough $(P_4 \cup P_{10})$-free graph on at least three vertices which has no 2-factor. 
Then we can take a special Tutte pair $(S,T)$ for $G$ by Theorem \ref{thm:2-factor}.
Note that $\delta(G) \geq 3$ since $G$ is a $\frac{3}{2}$-tough non-complete graph.

\begin{cla}\label{cla:size of S,T}
    $S \neq \emptyset$ and $|T| \geq 2$.
\end{cla}
\begin{prf}
	The expression $\eta(S,\emptyset)=2|S|$ immediately shows that $T$ is not empty. 
	Suppose that $S=\emptyset$.
	Then we have $h(\emptyset, T) \geq \delta(G) \geq 3$ by Lemma \ref{lem:sp Tutte pair1}$(ii)$. 
	Since $G$ is $\frac{3}{2}$-tough, $h(\emptyset,T) \leq \omega(G-T) \leq \frac{2}{3}|T|$ holds.
	However, since $(\emptyset,T)$ is a Tutte pair, we have
	\[-2 \geq \eta(\emptyset, T)=-2|T|+\sum_{x \in T} d_G(x)-h(\emptyset,T) \geq |T|-h(\emptyset,T) \geq \frac{|T|}{3} \geq \frac{1}{3},\]
	a contradiction. Thus, $S$ is not empty.

	Suppose that $|T|=1$, and let $T=\{x\}$. 
	Then by Lemma \ref{lem:sp Tutte pair1}$(ii)$, $d_{G-S}(x)=h(S,T)$ holds.
	However, since $S$ is not empty and $(S,T)$ is a Tutte pair, we have
	\[-2 \geq \eta(S,T) = 2|S|-2+d_{G-S}(x)-h(S,T)=2|S|-2 \geq 0,\]
	a contradiction. 
	Thus, we conclude $|T| \geq 2$. \qed
\end{prf}

In the following, the special Tutte pair $(S,T)$ is fixed.
Thus, we let $h=h(S,T)$.
Let 
\begin{eqnarray*}
    \begin{aligned}
        m&=\max\{d_{G-S}(x) \mid x \in T \},\text{ and }\\
        T_j&=\{x \in T \mid d_{G-S}(x)=j\}\text{ for } j=0,1,\ldots,m.
    \end{aligned}     
\end{eqnarray*}
Clearly, 
\[|T|=\sum_{j=0}^m|T_j| \text{ and } \sum_{x \in T}d_{G-S}(x)=\sum_{j=0}^m j|T_j|\]
hold.
Therefore, we have
\begin{equation}\label{eta}
    \eta(S,T)=2|S|+\sum_{j=0}^m(j-2)|T_j|-h \leq -2.
\end{equation}

Let $\mathcal{D}$ be the set of components of $G-S$ intersecting with $T \setminus T_0$ and $\mathcal{C}$ be the set of odd components of $G-(S \cup T)$. 
For each $D \in \mathcal{D}$ and $x \in T \setminus T_0$, we define 
\[\mathcal{C}_D=\{C \in \mathcal{C} \mid C \text{ is included in }D\} \text{ and }\mathcal{C}_x=\{C \in \mathcal{C} \mid e_G(x,C)=1\}.\]
Note that for any $x \in V(D) \cap T$ where $D \in \mathcal{D}$, the set $\mathcal{C}_x$ is a subset of $\mathcal{C}_D$.

Let $C \in \mathcal{C}$.
If $e_G(C,T) \geq 3$, we say that $C$ is {\it strong}; otherwise, $C$ is {\it weak}.
Note that if $C$ is weak, then $e_G(C,T)=1$.
For each $D \in \mathcal{D}$, we take a vertex $x_D \in V(D) \cap T$ such that the number of strong odd components adjacent to $x_D$ is maximum. 

\begin{case}\label{case1}
    For any $D \in \mathcal{D}$, every odd component in $\mathcal{C}_D \setminus \mathcal{C}_{x_D}$ is weak.
\end{case}

In this case, we define 
\begin{eqnarray*}
    \begin{aligned}
        T^\Rone &= \{x_D \mid D \in \mathcal{D}\},\\
        T^\Rone_j &= T^\Rone \cap T_j \text{ for }j=1,\ldots,m,\\
        T^\Rtwo &= \left\{y \in T  \mid y \text{ is adjacent to some }C \in \mathcal{C} \setminus \left(\bigcup_{D \in \mathcal{D}}\mathcal{C}_{x_D}\right)\right\}, \text{ and }\\
        T^\Rtwo_j &= T^\Rtwo \cap T_j \text{ for }j=1,\ldots,m.
    \end{aligned}     
\end{eqnarray*}

Let $D \in \mathcal{D}$ and $y \in T^\Rtwo \cap V(D)$.
Since $x_D$ and $y$ are in the same component $D$, we have $\mathcal{C}_{x_D} \cap \mathcal{C}_y \neq \emptyset$, which implies
\begin{equation}\label{h_upper in case1}
    h = \sum_{D \in \mathcal{D}} \left(|\mathcal{C}_{x_D}|+\sum_{y \in T^\Rtwo \cap V(D)}|\mathcal{C}_y \setminus \mathcal{C}_{x_D}| \right) \leq \sum_{j=1}^m \left( j |T^\Rone_j|+(j-1)|T^\Rtwo_j| \right).
\end{equation}

Let
\[W = S \cup T^\Rtwo \cup N_G \left( T_1 \setminus T^\Rone_1 \right).\]
Then each vertex in $T_0 \cup (T_1 \setminus T_1^\Rone)$ is an isolated vertex in the graph $G-W$.
Also, for any $D \in \mathcal{D}$, each component in $\mathcal{C}_D \setminus \mathcal{C}_{x_D}$ is separated from $D$ in the graph $G-W$ by the assumption of Case 1.
Therefore, since $|\mathcal{D}|=|T^\Rone|$ and $\sum_{D \in \mathcal{D}} |\mathcal{C}_{D} \setminus \mathcal{C}_{x_D}| \geq |T^\Rtwo|$, the graph $G-W$ has at least 
\[|T_0| + |T_1 \setminus T_1^\Rone| + |T^\Rone| + | T^\Rtwo | = |T_0| + |T_1| + |T^\Rone \setminus T_1| + | T^\Rtwo |\] 
components.
Moreover, (\ref{eta}) and (\ref{h_upper in case1}) imply
\begin{eqnarray}\label{size W in case1}
    \begin{aligned}
        |W| &= |S|+|T_1|+|T^\Rtwo| - |T^\Rone_1|\\
        & \leq |T_0|+\frac{3}{2}|T_1|+\frac{1}{2}\left(h-\sum_{j=2}^m(j-2)|T_j|\right)-1 + |T^\Rtwo| - |T^\Rone_1|\\
        & \leq |T_0|+\frac{3}{2}|T_1|+\frac{1}{2}\sum_{j=2}^m \left( 2|T^\Rone_j|+|T^\Rtwo_j|\right)-1 + |T^\Rtwo| - \frac{|T^\Rone_1|}{2}\\
        & \leq |T_0|+\frac{3}{2}|T_1|+|T^\Rone \setminus T_1|+\frac{3}{2}|T^\Rtwo|-\frac{| T^\Rone_1|}{2}-1\\
        &< |T_0|+\frac{3}{2}|T_1|+|T^\Rone \setminus T_1|+\frac{3}{2} |T^\Rtwo|.
    \end{aligned}
\end{eqnarray}

Suppose that $\omega(G-W)=1$.
Then since $|T| \geq 2$ by Claim \ref{cla:size of S,T}, $T_0 = T_1 = T^\Rtwo = \emptyset$ and $|T^\Rone|=1$ hold.
Let $T^\Rone = \{x\}$.
Then all odd components of $G-(S \cup T)$ are adjacent to $x$, which implies $m \leq h=d_{G-S}(x) \leq m$.
However, (\ref{eta}) follows
\[-2 \geq \eta(S,T)=2|S|+\sum_{j=2}^m(j-2)|T_j|-h \geq 2|S|+m-2-m=2|S|-2\]
implying $|S| = 0$, contrary to Claim \ref{cla:size of S,T}.
Therefore, $\omega(G-W) \geq 2$ holds.
Since the graph $G$ is $\frac{3}{2}$-tough, we have
\begin{eqnarray*}
    \begin{aligned}
        \frac{3}{2}\left( |T_0| + |T_1| + |T^\Rone \setminus T_1| + |T^\Rtwo| \right) \leq \frac{3}{2}\omega(G-W) \leq |W|,
    \end{aligned}
\end{eqnarray*}
a contradiction to (\ref{size W in case1}).

\begin{case}
    There exists $D \in \mathcal{D}$ such that the set $\mathcal{C}_D \setminus \mathcal{C}_{x_D}$ has some strong odd components.
\end{case}

We use the following terminology for paths.
Let $P$ be a path. 
If $u$ and $v$ are the end-vertices of $P$, we say that $P$ is a $uv$-path. 
For $x,y \in V(P)$, $xPy$ denotes the subpath of $P$ between $x$ and $y$. 
Let $P$ be an $xy$-path and $Q$ be a $yz$-path. 
If $P$ and $Q$ are internally disjoint, then $xPyQz$ denotes a path between $x$ and $z$ passing through $P$ and $Q$.
For a set $U$, if the end-vertices of a path $P$ are both in $U$, we say that $P$ is a $U$-path.  
The term defined below plays an important role in the proof.

\begin{defn}\label{def:basic path}
    Let $C \in \mathcal{C}$ be strong.
    For $U \subset N_T(C)$ with $|U| \geq 2$, a $U$-path $P$ is called a {\it basic $U$-path for $C$} if
    \begin{itemize}
        \item all internal vertices of $P$ are on $C$, and
        \item $P$ is an induced path in $G$ passing through at most two vertices in $N_C(U)$.
    \end{itemize} 
\end{defn}

In Definition \ref{def:basic path}, for any vertex $u \in U$, we can take a basic $U$-path for $C$ starting in $u$.
In fact, we can just choose the shortest path that starts from $u$ and goes through $C$ to another vertex of $U$.
Moreover, by Lemma \ref{lem:sp Tutte pair1}$(iii)$, every basic $U$-path for $C$ is on at least 4 vertices.

Let $D \in \mathcal{D}$ and $x=x_D$.
Suppose that $\mathcal{C}_D \setminus \mathcal{C}_x$ has some strong odd components.
Then we can take a strong odd component $C_0 \in \mathcal{C}_D \setminus \mathcal{C}_x$ such that there exists $C_1 \in \mathcal{C}_x$ with $N_T(C_0) \cap N_T(C_1) \neq \emptyset$.
For $U_1=\{x\} \cup (N_T(C_0) \cap N_T(C_1))$, let $Q_1$ be a basic $U_1$-path for $C_1$ starting in $x$ and $y \in N_T(C_0) \cap N_T(C_1)$ be the other end-vertex of $Q_1$.
Now we assume that the strong odd component $C_0 \in \mathcal{C}_D \setminus \mathcal{C}_x$ is chosen so that the number of vertices on $Q_1$ is as small as possible.

Since $C_0$ is strong, we can take a basic $N_T(C_0)$-path $Q_0$ for $C_0$ starting in $y$.
Let $z \in N_T(C_0) \setminus \{y\}$ be the other end-vertex of $Q_0$.
Also, by the choice of $x$ and $C_0 \in \mathcal{C}_y \setminus \mathcal{C}_x$, there exists a strong odd component $C_2 \in \mathcal{C}_x \setminus \mathcal{C}_y$. 
We take a basic $(N_T(C_2) \setminus \{z\})$-path $Q_2$ for $C_2$ starting in $x$ and let $w \in N_T(C_2) \setminus \{x,z\}$ be the other end-vertex of $Q_2$.
By the choice of those paths, Claim \ref{cla:taking} immediately follows.

\begin{cla}\label{cla:taking}
    Each of the following holds.
    \begin{enumerate}[$(i)$]
        \item Let $C \in \mathcal{C}_D \setminus \mathcal{C}_x$. 
        If $C$ is strong, then for any $u \in N_T(C) \setminus \{y\}$ we have $e_G(u,Q_1)=0$.
        \item For any $u \in T \setminus \{y,z\}$ we have $e_G(u,Q_0)=0$.
        \item For any $u \in T \setminus \{w,x,z\}$ we have $e_G(u,Q_2)=0$.
    \end{enumerate}
\end{cla}

\begin{figure}[t]
            \centering
            \includegraphics{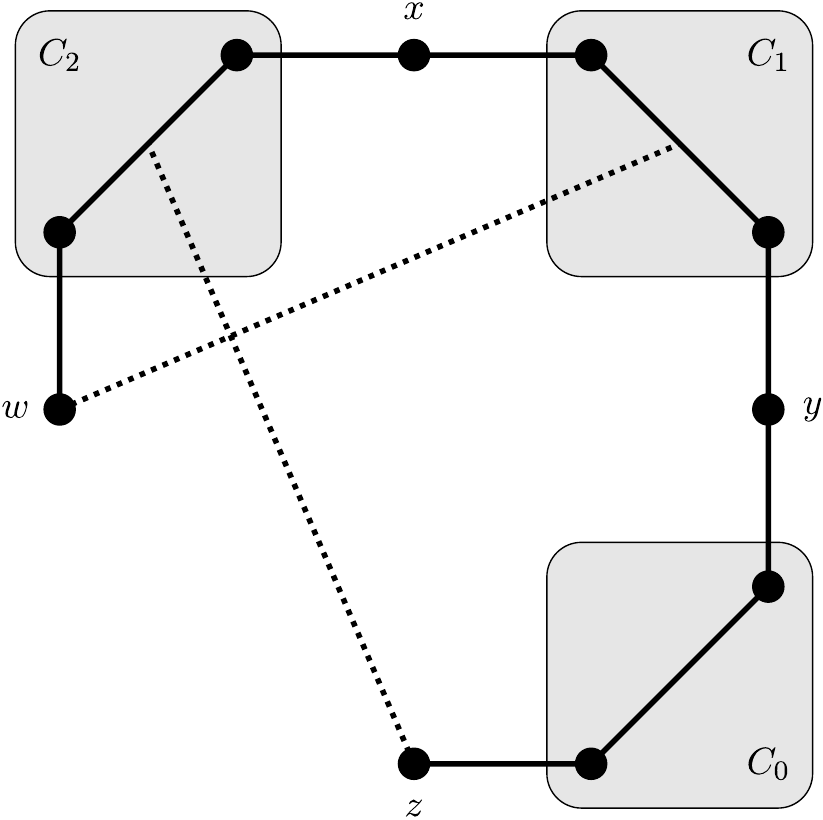}
        
            \caption{A sketch of the path $Q$. The dotted lines can be edges of $G$.}\label{fig:path_Q}
        \end{figure}

We define a path $Q$ in $G$ as
\begin{eqnarray*}
    Q= \left\{
        \begin{aligned}
            &w Q_2 x Q_1 y Q_0 z & \text{ if } &e_G(w,Q_1)=0 \text{ and } e_G(z,Q_2)=0;\\
            &w Q_2 z_2 z Q_0 y Q_1 x & \text{ if } &e_G(w,Q_1)=0 \text{ and } e_G(z,Q_2)=1;\\
            &x Q_2 w w_1 Q_1 y Q_0 z & \text{ if } &e_G(w,Q_1)=1 \text{ and } e_G(z,Q_2)=0;\\
            &x Q_1 w_1 w Q_2 z_2 z Q_0 y & \text{ if } &e_G(w,Q_1)=1 \text{ and } e_G(z,Q_2)=1,
        \end{aligned}
         \right.
\end{eqnarray*}
where 
\begin{eqnarray*}
    \begin{aligned}
        &z_2 \in N_G(z) \cap V(Q_2) &\text{ if }& e_G(z,Q_2)=1, \text{ and }\\
        &w_1 \in N_G(w) \cap V(Q_1) &\text{ if }& e_G(w,Q_1)=1.
    \end{aligned}     
\end{eqnarray*}
The path $Q$ is sketched in Figure \ref{fig:path_Q}. 
By Claim \ref{cla:taking}, $Q$ is an induced path in $G$ on at least 10 vertices.

\begin{cla}\label{cla:strong cpt}
    Let $C \in \mathcal{C}$ and $v \in T \setminus (T_0 \cup T_1)$.
    Then each of the following holds.
    \begin{enumerate}[$(i)$]
        \item If $C$ is strong, then $C$ is adjacent to at least one of $\{x,y,z,w\}$.
        \item If $v$ is adjacent to some weak odd component in $\mathcal{C}$, then $v$ is adjacent to at least one of $\{C_0,C_1,C_2\}$.
    \end{enumerate}
\end{cla} 
\begin{prf}
    For $(i)$, suppose that there exists a strong odd component $C \in \mathcal{C}$ such that $C$ is not adjacent to $x$, $y$, $z$ and $w$.
    Then we can take a basic $N_T(C)$-path $P$ for $C$, which is an induced path in $G$ on at least 4 vertices.
    Let $u,v \in N_T(C)$ be the end-vertices of $P$.
    By Claim \ref{cla:taking}, both $u$ and $v$ are not adjacent to any vertex on $Q$.
    Thus, $P \cup Q$ would contain $P_4 \cup P_{10}$ which is an induced subgraph of $G$, contrary to the assumption. 

    Next for $(ii)$, suppose that there exists a vertex $v \in T \setminus (T_0 \cup T_1)$ such that $v$ is adjacent to some weak odd component $C \in \mathcal{C}$, but not to $C_0$, $C_1$ and $C_2$.
    By Lemma \ref{lem:sp Tutte pair1}$(ii)$, we can take an odd component $C' \in \mathcal{C}_v \setminus \{C\}$.
    Also, by Lemma \ref{lem:sp Tutte pair1}$(iv)$, the odd component $C$ has at least 3 vertices.
    Let $p \in N_C(v)$, $q \in N_{C'}(v)$ and $r \in N_C(p)$.
    Then $P=r p v q$ is an induced path in $G$ on 4 vertices.
    Thus, $P \cup Q$ would contain $P_4 \cup P_{10}$ which is an induced subgraph of $G$, contrary to the assumption. \qed
\end{prf}

We define 
\begin{eqnarray*}
    \begin{aligned}
        \mathcal{C}^* &= \mathcal{C} \setminus (\mathcal{C}_x \cup \mathcal{C}_y \cup \mathcal {C}_z \cup \mathcal{C}_w),\\
        T^\Rone &= \{x,y,z,w\},\\
        T^\Rone_j &= T^\Rone \cap T_j \text{ for }j=1,\ldots,m,\\
        T^\Rtwo&=\{u \in T \mid \text{$u$ is adjacent to some }C \in \mathcal{C^*}\}, \text{ and }\\
        T_j^\Rtwo&=T^\Rtwo \cap T_j \text{ for }j=1,\ldots,m.
    \end{aligned}     
\end{eqnarray*}
The definitions of $\mathcal{C}^*$ and $T^\Rtwo$ imply $T^\Rone \subset T \setminus T^\Rtwo$. 
By Claim \ref{cla:strong cpt}$(i)$, every odd component in $\mathcal{C}^*$ is weak. 
Moreover, by Claim \ref{cla:strong cpt}$(ii)$, every vertex in $T^\Rtwo \setminus T_1^\Rtwo$ is adjacent to at least one of $\{C_0,C_1,C_2\}$.
Thus, we obtain 
\[|\mathcal{C}^*| \leq |T_1^\Rtwo|+\sum_{j=2}^m (j-1)|T_j^\Rtwo|,\] 
which implies
\begin{eqnarray}\label{h_case2}
    \begin{aligned}
        h &= |\mathcal{C}|\\ 
        &= |\mathcal{C}^* \cup \mathcal{C}_x \cup \mathcal{C}_y \cup \mathcal {C}_z \cup \mathcal{C}_w|\\
        &\leq |\mathcal{C}^*|+|\mathcal{C}_x|+|\mathcal{C}_y \setminus \{C_1\}|+|\mathcal{C}_z \setminus \{C_0\}|+|\mathcal{C}_w \setminus \{C_2\}|\\
        &\leq |T_1^\Rtwo|+\sum_{j=2}^m(j-1)|T_j^\Rtwo|+\sum_{j=1}^m j|T_j^\Rone|-3.
    \end{aligned}
\end{eqnarray}
Let 
\[W=S \cup (T^\Rtwo \setminus T_1^\Rtwo) \cup N_G(T_1 \setminus T_1^\Rtwo).\]
Then each vertex in $T_0 \cup (T_1 \setminus T_1^\Rtwo)$ is an isolated vertex in the graph $G-W$.
Also, each component in $\mathcal{C}^*$ is separated from the component $D$ in the graph $G-W$. 
Therefore, since $|\mathcal{C}^*| \geq |T^\Rtwo|$, the graph $G-W$ has at least 
\[|T_0| + |T_1 \setminus T_1^\Rtwo| + |T^\Rtwo| + 1 = |T_0|+|T_1|+|T^\Rtwo \setminus T_1^\Rtwo|+1\] 
components.
Moreover, (\ref{eta}) and (\ref{h_case2}) imply
\begin{eqnarray}\label{size W in case2}
    \begin{aligned}
        |W| &= |S|+|T^\Rtwo \setminus T^\Rtwo_1|+|T_1 \setminus T_1^\Rtwo|\\
        & \leq |T_0|+\frac{3}{2}|T_1|+\frac{1}{2}\left(h-\sum_{j=2}^m(j-2)|T_j|\right)-1 + |T^\Rtwo \setminus T^\Rtwo_1| - |T^\Rtwo_1|\\
        & \leq |T_0|+\frac{3}{2}|T_1|+\frac{1}{2}\left(\sum_{j=2}^m \left(2|T_j^\Rone|+|T_j^\Rtwo|\right)+|T_1^\Rone|-3\right) -1 + |T^\Rtwo \setminus T^\Rtwo_1| - \frac{|T^\Rtwo_1|}{2}\\
        & \leq |T_0|+\frac{3}{2}|T_1|+\frac{3}{2} + \frac{3}{2}|T^\Rtwo \setminus T^\Rtwo_1| - \frac{|T^\Rtwo_1|}{2}\\
        & \leq \frac{3}{2}(|T_0|+|T_1|+|T^\Rtwo \setminus T_1^\Rtwo|+1).
    \end{aligned}
\end{eqnarray}

\begin{subcase1}
    $\omega(G-W) \geq 2$.
\end{subcase1}

Since the graph $G$ is $\frac{3}{2}$-tough, we have
\begin{eqnarray*}\label{formula_case2}
    \begin{aligned}
        \frac{3}{2}(|T_0|+|T_1|+|T^\Rtwo \setminus T_1^\Rtwo|+1) &\leq \frac{3}{2}\omega(G-W) \leq |W|,\\
    \end{aligned}
\end{eqnarray*}
and so (\ref{h_case2}) and (\ref{size W in case2}) achieve the equations.
Therefore, we immediately see that the following statements hold.

\begin{itemize}
    \item $T_0$, $T_1^\Rone$ and $T_1^\Rtwo$ are empty.
    In particular, $z,w \notin T_1$, which implies that $\mathcal{C}_z \setminus \{C_0\}$ and $\mathcal{C}_w \setminus \{C_2\}$ are not empty.
    \item The sets $\mathcal{C}_x$, $\mathcal{C}_y \setminus \{C_1\}$, $\mathcal{C}_z \setminus \{C_0\}$ and $\mathcal{C}_w \setminus \{C_2\}$ are pairwise disjoint.
    Thus, we get $Q=w Q_2 x Q_1 y Q_0 z$.
\end{itemize}

Let $C_3 \in \mathcal{C}_z \setminus \{C_0\}$, $C_4 \in \mathcal{C}_w \setminus \{C_2\}$, $z_0 \in N_{C_0}(z)$ and $w_4 \in N_{C_4}(w)$.
Suppose that $C_3$ is strong. 
Then we can take a basic $(N_T(C_3) \setminus \{z\})$-path $Q_3$ for $C_3$.
Note that $Q_3$ is not adjacent to $Q_0$, $Q_1$ and $Q_2$ by Claim \ref{cla:taking}.
Thus, $Q_3 \cup (w_4 w Q z_0)$ would contain $P_4 \cup P_{10}$ which is an induced subgraph of $G$, contrary to the assumption.
Therefore, all odd components in $\mathcal{C}_z \setminus \{C_0\}$ are weak, and likewise all odd components in $\mathcal{C}_w \setminus \{C_2\}$ are weak.
These statements imply
\begin{equation}\label{contradict1}
    \omega(G-W-\{z,w\}) \geq |T_1|+|T^\Rtwo|+d_{G-S}(z)+d_{G-S}(w)-1 \geq 2.
\end{equation}
On the other hand, since $G$ is $\frac{3}{2}$-tough, we have
\begin{eqnarray}\label{contradict2}
    \begin{aligned}
        \frac{3}{2}\omega(G-W-\{z,w\}) \leq |W|+2 = \frac{3}{2}(|T_1|+|T^\Rtwo|+1)+2.
    \end{aligned}
\end{eqnarray}
However, (\ref{contradict1}) and (\ref{contradict2}) imply $d_{G-S}(z)+d_{G-S}(w) \leq \frac{10}{3} < 4$, a contradiction to $z,w \notin T_1$.

\begin{subcase2}
    $\omega(G-W)=1$.
\end{subcase2}

In this case, the sets $T_0$, $T_1$, $T^\Rtwo$ and $\mathcal{C}^*$ are empty.
Since $W=S$, Claim \ref{cla:size of S,T} and (\ref{size W in case2}) implies $|S|=1$. 
Therefore, (\ref{eta}) implies
\begin{eqnarray}\label{h_lower_bound}
    h \geq 4+\sum_{j=2}^m (j-2)|T_j| \geq \sum_{j=2}^m j|T_j^\Rone|-4.
\end{eqnarray}
Using the toughness assumption, we prove Claims \ref{cla:all odd cpts are strong} and \ref{cla:basic path}.

\begin{cla}\label{cla:all odd cpts are strong}
    All odd components of $G-(S \cup T)$ are strong.
\end{cla}
\begin{prf}
    Suppose that $C \in \mathcal{C}$ is weak.
    Then $C$ is adjacent to exactly one vertex $v \in T^\Rone$.
    Since the set $S \cup \{v\}$ separates $C$ from the others in $G$, the graph $G-S-v$ has at least two components.
    However, it implies
    \[\omega(G-S-v) \geq 2 > \frac{4}{3} = \frac{2}{3}(|S|+1),\]
    a contradiction to the toughness assumption. \qed
\end{prf}

\begin{cla}\label{cla:basic path}
    Let $C \in \mathcal{C}$. 
    Then for any vertex $v \in N_T(C)$ there exists a basic $N_T(C)$-path avoiding $N_C(v)$. 
\end{cla}
\begin{prf}
    Let $C \in \mathcal{C}$ and suppose that there is no basic $N_T(C)$-path avoiding $N_C(v)$ for a vertex $v \in N_T(C)$.
    Since $C$ is strong by Claim \ref{cla:all odd cpts are strong}, we can take two vertices $a,b \in N_T(C) \setminus \{v\}$.
    Let $u \in N_C(v)$.
    Then by the hypothesis, the vertices in $N_C(T)$ are separated in $C-u$.
    However, it implies 
    \[\omega(G-S-\{a,b,u\}) \geq 3 >\frac{8}{3}=\frac{2}{3}(|S|+3),\]
    a contradiction to the toughness assumption. \qed
\end{prf}

Using Claim \ref{cla:basic path} for some odd components, we will show that the graph $G$ contains $P_4 \cup P_{10}$ as an induced subgraph.
First, we show that the graph $G-(S \cup T)$ has an odd component which is adjacent to both $z$ and $w$.

\begin{cla}\label{cla:avoiding}
    For any $C \in \mathcal{C} \setminus \mathcal{C}_x$ we have $e_G(C, \{y,z,w\}) \geq 2$.
\end{cla}
\begin{prf}
    Let $C \in \mathcal{C} \setminus \mathcal{C}_x$.
    Then $C$ is adjacent to at least one of $\{y,z,w\}$ by Claim \ref{cla:strong cpt}.
    Suppose that $e_G(C, \{y,z,w\})=1$, and let $v \in \{y,z,w\}$ with $e_G(v,C)=1$.
    Then by Claim \ref{cla:basic path}, we can take a basic $N_T(C)$-path $P$ avoiding $N_C(v)$, which is an induced path of $G$ on at least 4 vertices.
    By Claim \ref{cla:taking}, we have $e_G(P,Q)=0$.
    Thus, $P \cup Q$ would contain $P_4 \cup P_{10}$ which is an induced subgraph of $G$, contrary to the assumption. \qed
\end{prf}

\begin{cla}\label{cla:an odd cpt with z and w}
    There exists an odd component $C_3 \in \mathcal{C} \setminus \{C_0,C_1,C_2\}$ which is adjacent to both $z$ and $w$.
\end{cla}
\begin{prf}
    Note that by (\ref{h_lower_bound}), there is at most one odd component of $G-(S \cup T)$ that belongs to at least two of the sets $\mathcal{C}_x$, $\mathcal{C}_y \setminus \{C_1\}$, $\mathcal{C}_z \setminus \{C_0\}$ and $\mathcal{C}_w \setminus \{C_2\}$.
    If $w$ is adjacent to some $C \in \{C_0,C_1\}$, then there exists $C' \in \mathcal{C}_z \setminus \{C_0,C_1,C_2\}$ with $e_G(C',\{x,y,w\})=0$, contrary to Claim \ref{cla:avoiding}.
    Thus, we get $e_G(w,C_0)=e_G(w, C_1)=0$, implying that $\mathcal{C}_w \setminus \{C_0,C_1,C_2\}$ is not empty.
    Let $C_3 \in \mathcal{C}_w \setminus \{C_0,C_1,C_2\}$.
    Then one of $\{x,y,z\}$ is adjacent to $C_3$ by Claim \ref{cla:avoiding}. 
    If $e_G(C_3,\{x,y\}) \geq 1$, then there exists $C \in \mathcal{C}_z \setminus \{C_0,C_1,C_2,C_3\}$ with $e_G(C,\{x,y,w\})=0$, contrary to Claim \ref{cla:avoiding}.
    Thus, we obtain $e_G(C_3,\{x,y\})=0$ implying $e_G(C_3,z)=1$.\qed
\end{prf}

By Claim \ref{cla:an odd cpt with z and w}, we can take an odd component $C_3 \in (\mathcal{C}_z \cap \mathcal{C}_w) \setminus \{C_0,C_1,C_2\}$.
Thus, we obtain 
\begin{eqnarray*}
    \begin{aligned}
        h&=|\mathcal{C}_x|+|\mathcal{C}_y \setminus \{C_1\}|+|\mathcal{C}_z \setminus \{C_0\}|+|\mathcal{C}_w \setminus \{C_2,C_3\}| \\
        &=\sum_{j=2}^m (j-2)|T_j|+4 = \sum_{j=2}^m j|T^\Rone_j|-4,
    \end{aligned}
\end{eqnarray*}
which implies that the following statements hold.

\begin{itemize}
    \item The sets $\mathcal{C}_x$, $\mathcal{C}_y \setminus \{C_1\}$, $\mathcal{C}_z \setminus \{C_0\}$ and $\mathcal{C}_w \setminus \{C_2,C_3\}$ are pairwise disjoint. Thus, we get $Q=w Q_2 x Q_1 y Q_0 z$.
    \item $T \setminus T^\Rone \subset T_2$.
\end{itemize}

We shall show that the set $T \cup V(C_0) \cup V(C_1) \cup V(C_2) \cup V(C_3)$ contains a subset which induces a subgraph of $G$ isomorphic to $P_4 \cup P_{10}$.

Using Claim \ref{cla:basic path}, we can take a basic $N_T(C_3)$-path $R_3$ avoiding $N_{C_3}(z)$.
Let $a,b \in N_{T}(C_3) \setminus \{z\}$ be the end-vertices of $R_3$.
Note that any vertex in $\{a,b\} \setminus \{w\}$ is not adjacent to $Q$ by Claim \ref{cla:taking}.
Since $G$ is $(P_4 \cup P_{10})$-free, $w \in \{a,b\}$ and we may set $b=w$.
Let $Q'=a R_3 w Q y$.
Then $Q'$ is an induced path in $G$ on at least 10 vertices.

Using Claim \ref{cla:basic path} again, we can take a basic $N_T(C_0)$-path $R_0$ avoiding $N_{C_0}(y)$.
Let $b,c \in N_{T}(C_0) \setminus \{y\}$ be the end-vertices of $R_0$ (note that the vertex $z$ may be in $\{b,c\}$).
Now any vertex in $\{b,c\} \setminus \{a\}$ is not adjacent to $Q'$.
Since $G$ is $(P_4 \cup P_{10})$-free, $a \in \{b,c\}$ and we may set $c=a$.
Replace $Q'$ with $b R_0 a Q' x$.
Then $Q'$ is still an induced path in $G$ on at least 10 vertices.

Using Claim \ref{cla:basic path} again, we can take a basic $N_T(C_1)$-path $R_1$ avoiding $N_{C_1}(x)$.
Let $c,d \in N_{T}(C_1) \setminus \{x\}$ be the end-vertices of $R_1$ (note that the vertex $y$ may be in $\{c,d\}$).
Since $a \in T \setminus T^\Rone \subset T_2$, $a \notin \{c,d\}$.
Moreover, any vertex in $\{c,d\} \setminus \{b\}$ is not adjacent to $Q'$.
Thus, since $G$ is $(P_4 \cup P_{10})$-free, $b \in \{c,d\}$ and we may set $d=b$ (also $b \neq z$ turns out).
Replace $Q'$ with $c R_1 b Q' w$.
Then $Q'$ is still an induced path in $G$ on at least 10 vertices.

Using Claim \ref{cla:basic path} again, we can take a basic $N_T(C_2)$-path $R_2$ avoiding $N_{C_2}(w)$.
Let $d,e \in N_{T}(C_2) \setminus \{w\}$ be the end-vertices of $R_2$ (note that the vertex $x$ may be in $\{d,e\}$).
Since $a,b \in T \setminus T^\Rone \subset T_2$, $\{d,e\} \cap \{a,b\}=\emptyset$.
Moreover, any vertex in $\{d,e\} \setminus \{c\}$ is not adjacent to $Q'$.
Thus, since $G$ is $(P_4 \cup P_{10})$-free, either $c \in \{d,e\}$ and we may set $e=c$ (also $c \neq y$ turns out).
Replace $Q'$ with $d R_2 c Q' a$.
Then $Q'$ is still an induced path in $G$ on at least 10 vertices.

\begin{figure}[!b]
        \centering
    \includegraphics{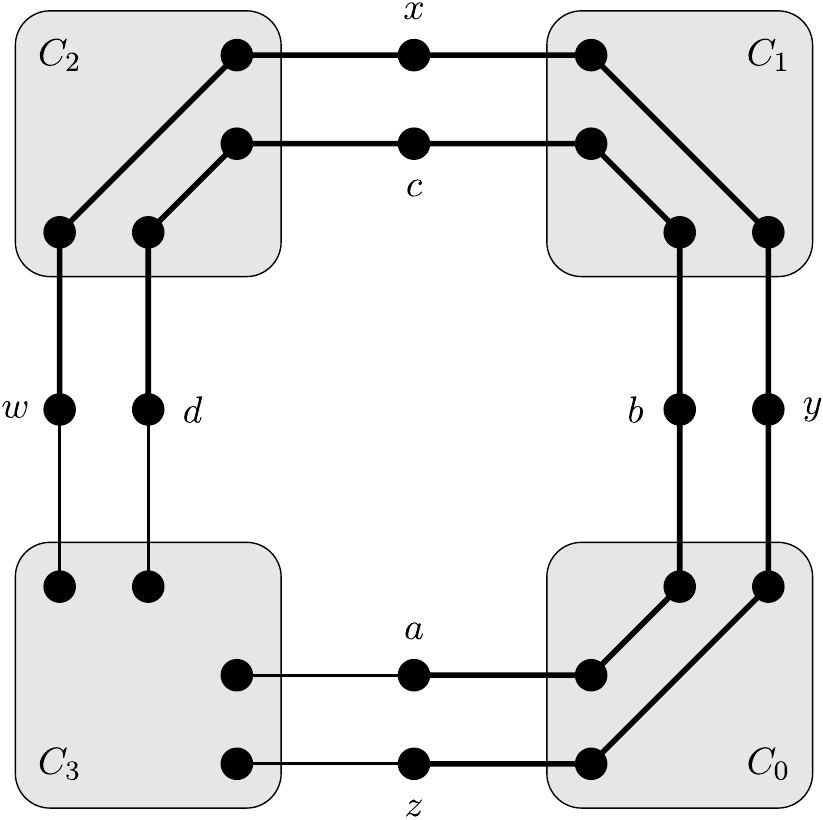}
    
            \caption{The situation around the paths $Q$ and $Q'$.}\label{fig:path_Q_Q'}
        \end{figure}

Suppose $d \notin N_T(C_3)$.
Using Claim \ref{cla:basic path}, we take a basic $N_T(C_3)$-path $R$ avoiding $N_{C_3}(a)$.
Let $e,f \in N_{T}(C_3) \setminus \{a\}$ be the end-vertices of $R$.
Then we have $d \notin \{e,f\}$ by the hypothesis.
However, then $R \cup Q'$ would contain $P_4 \cup P_{10}$ which is an induced subgraph of $G$, contrary to the assumption.
Thus, we obtain $d \in N_T(C_3)$ (also $d \neq x$ turns out). 
(The situation around the paths $Q$ and $Q'$ is shown in Figure \ref{fig:path_Q_Q'}.)

Since $e_G(C_3,T) \equiv 1 \mod 2$, there exists a vertex $e \in N_T(C_3) \setminus \{a,d,w,z\}$.
Then $e$ is adjacent to neither $Q$ nor $Q'$.
Let $R$ be a basic $N_T(C_3)$-path starting in $e$ and $f \in N_T(C_3)$ be the other end-vertex of $R$.
If $f \notin \{a,d\}$, then $R \cup Q'$ would contain $P_4 \cup P_{10}$ which is an induced subgraph of $G$. Otherwise, $R \cup Q$ would contain $P_4 \cup P_{10}$ which is an induced subgraph of $G$.
Both cases contradict the assumption.
This completes the proof of Theorem \ref{thm:main1}. \qed

\section{Examples of graphs without 2-factor}

In this section, we show the examples of $2P_5$-free graphs without a 2-factor mentioned in Section \ref{sec:intro}.
The examples appear in \cite{Bauer_et_al-1994}, but we will describe the details again here.

We denote by $G+H$ the join of two graphs $G$ and $H$, which is defined by $V(G+H)=V(G) \cup V(H)$ and $E(G+H)=E(G) \cup E(H) \cup \{uv \mid u \in V(G) \text{ and }v \in V(H)\}$.
For two integers $l \geq 1$ and $m \geq 2$, we construct the graph $G(l,m)$ as follows. 
Let $A_1,\ldots,A_{2m+1}$ be $2m+1$ copies of $K_{2l+1}$, and $B$ be a copy of $K_{(2l+1)(2m+1)}$.
Let $A$ be the union of $A_1,\ldots,A_{2m+1}$.
The graph $H(l,m)$ is obtained by first joining the vertices of $A$ to the vertices of $B$ by a perfect matching, and then subdividing each matching edge with a single vertex.
We write the set of $(2l+1)(2m+1)$ subdividing vertices by $T$.
We define the graph $G(l,m)$ as $G(l,m)=S+H(l,m)$, where $S$ is a copy of $K_m$.
The graph $G(l,m)$ is sketched in Figure \ref{fig:G(l,m)}.

For the pair of sets $(S,T)$, since every vertex $v \in T$ satisfies $d_{G(l,m)-S}(v)=2$ and $A_1,\ldots,A_{2m+1}$ and $B$ are odd components of $G(l,m)-(S \cup T)$, we have
\begin{eqnarray*}
    \begin{aligned}
        \eta(S,T)=2|S|-2|T|+\sum_{v \in T}d_{G(l,m)-S}(v)-h(S,T)=2m-(2m+2)=-2.
    \end{aligned}
\end{eqnarray*}
Thus, by Theorem \ref{thm:2-factor}$(ii)$, the graph $G(l,m)$ does not contain 2-factor for any $l$ and $m$.
Moreover, for any $l$ and $m$, the graph $G(l,m)$ is a $2P_5$-free graph.
Indeed, any induced path in $G(l,m)$ on at least 3 vertices is in $H(l,m)$ and any induced path in $H(l,m)$ on at least 5 vertices must contain at least one vertex in $B$.
Therefore, since $B$ is complete, the graph $G(l,m)$ does not contain induced subgraph isomorphic to $2P_5$.

    \begin{figure}[!b]
        \centering
        \includegraphics{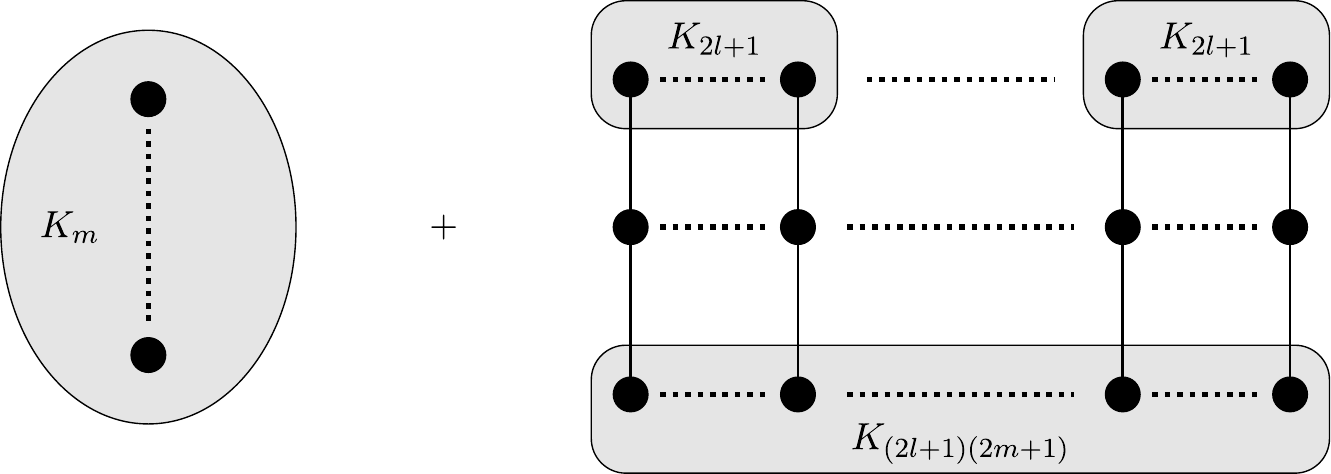}
    
        \caption{The graph $G(l,m)$.}\label{fig:G(l,m)}
    \end{figure}

The toughness of $G(l,m)$ is given by 
\[t(G(l,m))=\frac{m+(2l+1)(2m+1)+2l(2m+1)}{(2l+1)(2m+1)+1}=2-\frac{m+3}{(2l+1)(2m+1)+1}.\]
In particular, a set $W \subset V(G(l,m))$ with $\omega(G(l,m)-W) \geq 2$ that satisfies
\[t(G(l,m))=\frac{|W|}{\omega(G(l,m)-W)}\]
is obtained as follows.
First, for each $i=1,\ldots,2m+1$, we take one vertex $a_i \in A_i$.
Next, we choose one vertex $a \in \{a_1,\ldots,a_{2m+1}\}$, and let $x \in T$ be the vertex adjacent to $a$ and $b \in B$ be the vertex adjacent to $x$.
Then the set $W$ is obtained by
\[W=S \cup (A \setminus \{a_1,\ldots,a_{2m+1}\}) \cup (B \setminus \{b\}) \cup \{x\}.\]
By fixing $m$ and choosing $l$ large, we can make the toughness of the graph $G(l,m)$ approach to 2 from below.

\section*{Acknowledgements}
We would like to thank the referee and Professor Katsuhiro Ota for their helpful comments.
The author's work was supported by JST Doctoral Program Student Support Project (JPMJSP2123).

\bibliography{article2_12} 
\bibliographystyle{myplain} 

\end{document}